\documentclass{article}
\usepackage{amsmath,amsfonts,amsthm,epsfig,amssymb}
\newtheorem{theorem}{Theorem}

\newtheorem{corollary}[theorem]{Corollary}

\renewcommand\Pr{{\mathop{\mathbb P{}}\nolimits}}
\renewcommand\phi{\varphi}

\newcommand\Z{{\mathbb Z}}
\newcommand\RR{{\mathbb R}}
\newcommand\Gd{G^{\boxtimes}}
\newcommand\Gs{G^{\star}}
\newcommand\es{e^{\star}}
\newcommand\eps{\varepsilon}
\newcommand\cc[1]{{\overline{#1}}}
\newcommand\pH{p_{\mathrm H}}
\newcommand\pHb{p_{\mathrm H}^{\mathrm b}}
\newcommand\pcb{p_{\mathrm c}^{\mathrm b}}
\newcommand\bond{{\mathrm b}}
\newcommand\site{{\mathrm s}}
\newcommand\pHs{p_{\mathrm H}^{\mathrm s}}
\newcommand\bb[1]{\bigl(#1\bigr)}
\newcommand\cb[1]{\Bigl(#1\Bigr)}
\newcommand\lrb[1]{\left(#1\right)}

\begin{document}
\title{Percolation on dual lattices with $k$-fold symmetry}
\date{February 5, 2007}

\author{B\'ela Bollob\'as\thanks{Department of Pure Mathematics
and Mathematical Statistics, University of Cambridge, Cambridge CB3 0WB, UK}
\thanks{University of Memphis, Memphis TN 38152, USA}
\thanks{Research supported in part by NSF grants CCR-0225610,
 DMS-0505550 and W911NF-06-1-0076}
\and Oliver Riordan${}^*$%
\thanks{Royal Society Research Fellow}}
\maketitle

\begin{abstract}
Zhang found a simple, elegant argument deducing the non-existence
of an infinite open cluster in certain lattice percolation models (for example,
$p=1/2$ bond percolation on the square lattice) from general
results on the uniqueness of an infinite open cluster when it exists;
this argument requires some symmetry. Here we show that a simple modification
of Zhang's argument requires only $2$-fold (or $3$-fold) symmetry,
proving that the critical probabilities for percolation on dual planar
lattices with such symmetry sum to $1$. Like Zhang's argument,
our extension applies in many contexts; in particular,
it enables us to answer a question of Grimmett concerning
the anisotropic random cluster model on the triangular lattice.
\end{abstract}

\section{Introduction and results}

Let $G$ be an infinite, locally finite graph, and $0<p<1$ a real parameter.
In {\em independent bond percolation} on $G$, each edge of $G$ is assigned
a {\em state}, {\em open} or {\em closed}; the states of the edges are independent,
and each edge is open with probability $p$. An {\em open cluster} is a maximal
connected subgraph of $G$ all of whose edges are open, i.e., a component
of the {\em open subgraph} formed from $G$ by deleting the closed edges. We
write $\Pr_p$ for the corresponding probability measure.

Writing $\theta(p)=\theta_G(p)=\theta_{G;v}^\bond(p)$ for the probability that a fixed `starting vertex'
$v$ is in an infinite open cluster, 
the (Hammersley) {\em critical probability} $\pH=\pHb(G)$ 
for bond percolation on $G$ is defined by
\[
 \pHb(G) = \inf\{p: \theta_G(p)>0\},
\]
where $v$ is any fixed vertex of $G$.
(There are other natural notions of critical probability. For
well behaved graphs, including the lattices studied here, these
coincide, and one writes $\pcb$ for their common value. Here
we are concerned only with $\pHb$.)
Let $E_\infty$ denote the event that there is an infinite open cluster somewhere in $G$.
Kolmogorov's $0$/$1$-law implies that $\Pr_p(E_\infty)$ is either $0$ or $1$ for any $p$.
Simple standard arguments show that $\Pr_p(E_\infty)=0$ for $p<\pH$, and $\Pr_p(E_\infty)=1$
for $p>\pH$.

\begin{figure}[h]
\centering
\epsfig{file=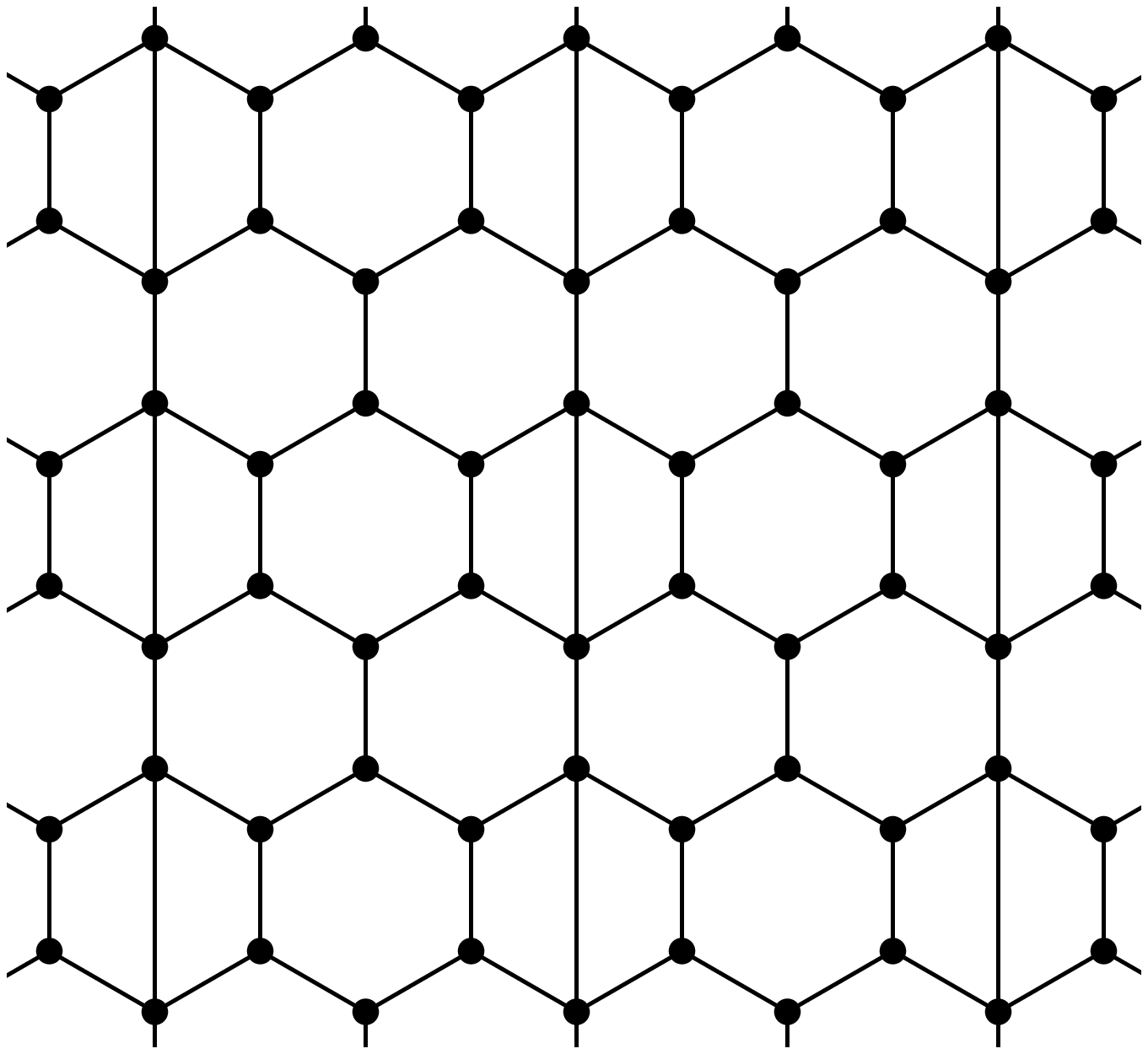,width=2in}
\hskip.5in
\epsfig{file=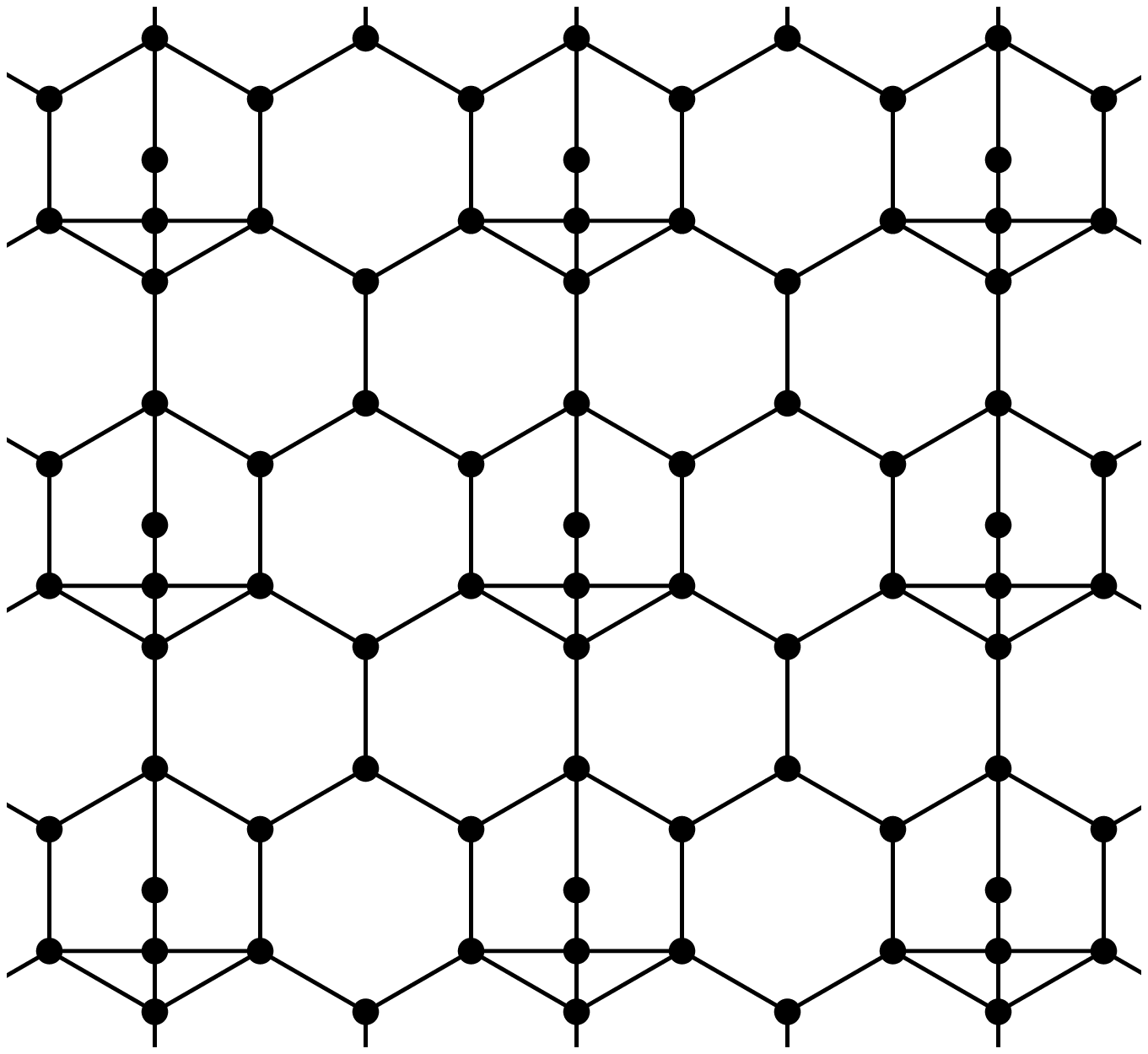,width=2in}
\caption{Two plane lattices obtained by `decorating' the hexagonal lattice.
The lattice on the left has $2$-fold symmetry if the origin is chosen suitably.
The lattice on the right has mirror symmetry, but does not have $k$-fold
symmetry for any $k\ge 2$.}\label{fig1}
\end{figure}

By a {\em plane lattice} we shall mean a (multi-)graph $G$ drawn in the plane
with non-crossing edges, with
the following properties: $G$ is connected, infinite, and locally finite,
the vertex set $V(G)$ is a subset of $\RR^2$ containing no accumulation points,
and there are independent vectors $w_1$, $w_2\in \RR^2$ such that the translations
$T_{w_i}$ of $\RR^2$ through the vectors $w_i$ induce isomorphisms of $G$ as a plane
graph. Usually, $G$ is drawn with straight edges, as in the examples in Figure~\ref{fig1},
but we shall not require this. However, we shall assume, as we may, that the
edges of $G$ are drawn as piecewise linear curves.

Note that the vertex set of a plane lattice $G$ need not be a lattice in the algebraic sense.
However, for every vertex $v$ of $G$ there is an isomorphism of $G$ mapping
$v$ to a vertex $v'$ in the parallelogram with corners $0$, $w_1$, $w_1+w_2$, $w_2$,
so $V(G)$ is the union of finitely many translates of the lattice $\{a_1w_1+a_2w_2:
a_1,a_2\in \Z\}$.

For $k\ge 2$,
we say that a plane lattice $G$ has {\em $k$-fold symmetry} if the rotation
about the origin through an angle of $2\pi/k$ maps
the plane graph $G$ into itself.
For example, the square lattice has $4$-fold symmetry. The hexagonal lattice
has $6$-fold symmetry as long as the origin is chosen appropriately,
i.e., as the centre of a hexagon. Figure~\ref{fig1} shows two plane lattices; the
one on the left has $2$-fold symmetry, the one on the right does not have
$k$-fold symmetry for any $k\ge 2$.

The {\em dual} $\Gs$ of a plane (multi-)graph $G$ is the (multi-)graph with one vertex
for every face of $G$, and an edge $\es$ for every edge of $G$;
the edge $\es$ joins the vertices of $\Gs$ corresponding to the faces
of $G$ in whose boundary $e$ lies. It is easy to see that,
if $G$ is a plane lattice, then so is an appropriate drawing of $\Gs$.
If, in addition, $G$ has $k$-fold symmetry, then $\Gs$ may be drawn with
$k$-fold symmetry.

Aizenman, Kesten and Newman~\cite{AKNunique} showed that, under rather
mild conditions, the probability that there is more than one infinite
open cluster is zero; a little later, Burton and Keane~\cite{BK} gave a very simple
proof of this result. This result applies in particular to independent
bond percolation on a plane lattice (as well as to lattices in all
dimensions). In 1988, Zhang found a simple way to deduce that if $G$ is a plane
lattice with $k$-fold symmetry, $k\ge 4$, and $0<p<1$, then one cannot have both
$\theta_G(p)>0$ and $\theta_{\Gs}(1-p)>0$; this argument
was first published in Grimmett~\cite[p.~195]{GG89}.
Zhang's argument
uses only symmetry, Harris's positive correlation
lemma, and the Aizenman--Kesten--Newman result.  For the square
lattice, which is self-dual, it follows that $\theta(1/2)=0$, giving an alternative proof
of this old result of Harris~\cite{Harris}. Using, for example, the exponential
decay result of Menshikov~\cite{Menshikov}, it follows that $\pH(\Z^2)=1/2$;
see \cite{GG99} or \cite[Section 5.2]{BRbook}.
The advantage of Zhang's argument over the (by now) many other
proofs of this celebrated result of Kesten~\cite{Kesten1/2}
(see~\cite{GG99,ourKesten,BRbook}, for example)
is that it does not rely on a `Russo--Seymour--Welsh'
type lemma~\cite{Russo,SW}, and so
adapts easily to other lattices, and also to other models,
such as the random cluster model of Fortuin and Kasteleyn~\cite{FK}.

Our aim is to show that, contrary to the generally held belief
(see Grimmett \cite{GG99},
for example), it is fairly easy to adapt Zhang's argument
to weaken the symmetry requirement: $k$-fold symmetry for any $k\ge 2$ is enough.

\begin{theorem}\label{th1}
Let $G$ be a plane lattice with $k$-fold symmetry for some $k\ge 2$,
and let $0<p<1$. Then either $\theta_G(p)=0$ or $\theta_{\Gs}(1-p)=0$.
\end{theorem}

We prove this result in the next section;
in the final section we shall discuss some extensions and an application.

Before turning to the proof, let us remark that,
in a paper published after the first draft of this paper was written,
Sheffield~\cite[Theorem 9.3.1]{Sheffield}
has proved a result corresponding to Theorem~\ref{th1}, but applying
to all plane lattices, with no additional symmetry requirement.
His proof of this stronger result is rather involved.

\section{The proof}

As usual in percolation, we declare an edge $\es$ of $\Gs$
to be open if and only if $e$ is closed. Note that the states of the edges
of $\Gs$ are independent, and each edge is open with probability $1-p$.

By an {\em open path} in $G$ we mean a finite or infinite
path $P=v_0e_0v_1e_1v_2\cdots $ in the graph $G$ such that
every edge $e_i$ of $P$ is open. An open path in $\Gs$ is
defined similarly. Since each edge of $G$ (or $\Gs)$ is a piecewise
linear curve in the plane, so is any open path $P$.

Let $D_R$ be the closed disc of radius $R$ centred at the origin.
We say that an infinite open path $P$ in $G$ or $\Gs$ {\em leaves $D_R$
at the point $x$} if, when $P$ is viewed as a curve in the plane,
the last point of $P$ in $D_R$ is $x$,
i.e., if $P$ contains an infinite piecewise linear path meeting $D_R$ only at $x$.
If the boundary of $D_R$ is divided into arcs $A_1,\ldots,A_r$,
then we write $L_i=L_i(R;A_i)$ for the event that there is an infinite open path
in $G$ leaving $D_R$ at some point $x\in A_i$.
Similarly, we write $L_i^\star$ for the event that there is an infinite
open path in $\Gs$ leaving $D_R$ at some point of $A_i$.
Note that each $L_i$ is an {\em increasing} event: if $L_i$ holds
and we change the states of one or more edges from closed to open, then $L_i$ still holds.

The essence of Zhang's argument is as follows. Let $G$ be a plane
lattice with $k$-fold symmetry, $k\ge 4$,
and suppose that $\theta_G(p)>0$ and $\theta_{\Gs}(1-p)>0$.
From Kolmogorov's $0$/$1$-law, with probability $1$ both $G$ and $\Gs$ contain
infinite open clusters.
As $R\to\infty$, the discs $D_R$ increase to cover $\RR^2$, so the probability
that $D_R$ meets an infinite open cluster in $G$ tends to $1$,
i.e., for any $\eps>0$, if $R$ is large enough this event has probability
at least $1-\eps$. Divide the boundary of $D_R$ symmetrically into
$k$ arcs $A_1,\ldots,A_k$. If $D_R$ meets an infinite open cluster in $G$,
then there is an infinite open path in $G$ leaving $D_R$ at some point $x$,
which must lie in one of the $A_i$. Thus, the union $L_1\cup\cdots\cup L_k$
has probability at least $1-\eps$.
Since the $L_i$ are increasing, it follows from Harris's Lemma (see below) that
$\Pr(L_i)\ge 1-\eps^{1/k}$ for some $i$. But $\Pr(L_i)=\Pr(L_j)$
by symmetry. Thus, if $R$ is large enough, we have $\Pr(L_i)\ge 1-\eps^{1/k}$
for every $i$.
Similarly, if $R$ is large enough, $\Pr(L_i^\star)\ge 1-\eps^{1/k}$ for every $i$.

Taking $\eps<(1/4)^k$, we find that with positive probability the event
$L_1\cap L_2^\star\cap L_3\cap L_4^\star$ holds. But simple topology shows that
when this event holds, one of $G$ and $\Gs$ contains two infinite open clusters.
This contradicts the Aizenman--Kesten--Newman Theorem.

To adapt this argument to the case $k=2$ or $k=3$, we simply divide the boundary
of $D_R$ into $2k$ arcs, in a way that has $k$-fold symmetry. Roughly speaking,
by moving the division point between $A_1$ and $A_2$ appropriately, we can
ensure that $L_1$ and $L_2$ (and hence, by symmetry, all $L_i$) have
about the same probability, so $\Pr(L_1)$ is close to $1$. With the arcs now
fixed, there is {\em some} $i$ with $\Pr(L_i^\star)$ close to $1$,
and hence $\Pr(L_{i+2}^\star)=\Pr(L_i^\star)$ also close to $1$,
so $L_i^\star\cap L_{i+1}\cap L_{i+2}^\star\cap L_{i+3}$ has positive probability.

For a reader familiar with Zhang's argument, it is a
straightforward exercise to turn the above outline into a complete proof of Theorem~\ref{th1};
nevertheless, we shall spell out the details in full.

\begin{proof}[Proof of Theorem~\ref{th1}]
Let $G$ be a plane lattice with $k$-fold symmetry, and let $\Gs$ be its dual,
drawn as a plane lattice with $k$-fold symmetry. Let $0<p<1$ be fixed,
and suppose for a contradiction that $\theta_G(p)>0$ and $\theta_{\Gs}(1-p)>0$.
We assume for convenience that $G$ and $\Gs$ are drawn with all edges
piecewise linear.
Rotating the coordinate system, if necessary, we may assume that no line segment
making up a part of an edge of $G$ or $\Gs$ is parallel to the $x$-axis.
Thus, the set $B_1$ of $R>0$ such that the point $(R,0)$ lies on an edge of $G$ or $\Gs$ is countable.
As $V(G)$ is countable, the set $B_2$ of $R>0$ such that any vertex of $G$
or $\Gs$ lies on the boundary of $D_R$ is also countable.
From now on we only consider $R\notin B=B_1\cup B_2$.

Given $R\notin B$ and $0\le\phi\le 2\pi/k$ such that $Q_1=(R\cos\phi,R\sin\phi)$ does not
lie on an edge of $G$ or $\Gs$, let $A_1$ be the boundary arc of $D_R$ from
$Q_0=(R,0)$ to $Q_1$, and let $A_2$ be the boundary arc of $D_R$ from $Q_1$
to $Q_2=(R\cos(2\pi/k),R\sin(2\pi/k))$.
For $i\ge 3$, let $A_i$ and $Q_i$ be defined by rotating $A_{i-2}$ and $Q_{i-2}$
anticlockwise about the origin through an angle $2\pi/k$; when referring to $A_i$ or $Q_i$,
we shall always take the subscript modulo $2k$.
Since no $Q_i$ lies on an edge of $G$ or $\Gs$, it will be irrelevant whether $A_i$
includes its endpoints; if we include one endpoint into each $A_i$, then the $2k$
arcs $A_i$ partition the boundary of $D_R$.

Let $L_i=L_i(\phi)=L_i(R,\phi)$ be the event that $G$ contains an infinite open path leaving $D_R$
at a point of $A_i$, and let $L_i^\star=L_i^\star(\phi)=L_i^\star(R,\phi)$ be defined similarly for $\Gs$.
Again, we take the subscript modulo $2k$.
Note that $\bigcup_i L_i$ is the event that $D_R$ meets
an infinite open cluster in $G$.

Since $\theta_G(p)>0$, the probability that $G$ contains an infinite open cluster
is positive, and hence (by Kolmogorov's $0$/$1$-law, for example) equal to $1$.
Thus, $\Pr(\bigcup_i L_i)\to 1$ as $R\to \infty$.
Set
\begin{equation}\label{edef}
 \eps=(1-p)^{2k}/5^{2k}>0,
\end{equation}
noting that $\eps$ depends on $p$
but not on $R$ or $\phi$. 
Then, if $R$ is large
enough, we have
\begin{equation}\label{uli}
 \Pr\cb{\bigcup_i L_i} \ge 1-\eps.
\end{equation}
Similarly, since $\theta_{\Gs}(1-p)>0$, if $R$ is large enough we have
\begin{equation}\label{ulis}
 \Pr\cb{\bigcup_i L_i^\star} \ge 1-\eps.
\end{equation}
For the rest of the proof, we fix an $R\notin B$ such that \eqref{uli} and \eqref{ulis} hold.

The events $L_i$ are increasing, so their complements $\cc{L_i}$ are decreasing.
Harris~\cite{Harris} showed that, if $E_1$ and $E_2$ are decreasing events in a product probability
space, then $\Pr(E_1\cap E_2)\ge \Pr(E_1)\Pr(E_2)$. As the intersection of two decreasing
events is decreasing, it follows that if $E_1,\ldots,E_j$ are decreasing,
then $\Pr(\bigcap_i E_i)\ge \prod_i\Pr(E_i)$. In particular,
\[
 \Pr\lrb{\bigcap_{i=1}^{2k} \cc{L_i}} \ge \prod_{i=1}^{2k} \Pr\bb{\cc{L_i}}
  = \Pr\bb{\cc{L_1}}^k\Pr\bb{\cc{L_2}}^k,
\]
where the equality follows from $k$-fold symmetry. From \eqref{uli}, it then
follows that
\begin{equation}\label{l12}
 \Pr\bb{\cc{L_1}}\Pr\bb{\cc{L_2}} \le \eps^{1/k},
\end{equation}
for any $0\le \phi\le 2\pi/k$.
Arguing in exactly the same way but for $\Gs$, it follows from \eqref{ulis}
that
\begin{equation}\label{l12s}
 \Pr\bb{\cc{L_1^\star}}\Pr\bb{\cc{L_2^\star}} \le \eps^{1/k}.
\end{equation}

As $G$ is drawn with piecewise linear edges, the set of boundary points
of $D_R$ that lie on edges of $G$ is countable. Furthermore,
as $G$ is a plane lattice, this set contains no accumulation points,
and is thus finite.
Let $0\le\phi<\phi'\le2\pi /k$ be such that there is exactly one point $x$ of the
arc of $D_R$ from $(R\cos\phi,R\sin\phi)$ to $(R\cos\phi',R\sin\phi')$
that lies on an edge of $G$, and suppose that $x$ lies in the interior of this arc.
Let $L_x$ be the event that $G$ contains an infinite open path leaving $D_R$
at the point $x$.
Then $L_1(\phi')=L_1(\phi)\cup L_x$, and $L_2(\phi)=L_2(\phi')\cup L_x$.
The events $\cc{L_i(\cdot)}$ and $\cc{L_x}$ are decreasing. Thus, from
Harris's Lemma,
\[
 \Pr\bb{\cc{L_1(\phi')}} = \Pr\bb{\cc{L_1(\phi)}\cap \cc{L_x}}
 \ge \Pr\bb{\cc{L_1(\phi)}} \Pr\bb{\cc{L_x}}.
\]
If the unique edge of $G$ on which $x$ lies is closed, then $\cc{L_x}$ holds,
so $\Pr\bb{\cc{L_x}}\ge 1-p$. Since $L_1(\phi)\subset L_1(\phi')$, we thus have
\begin{equation}\label{zz1}
 \Pr\bb{\cc{L_1(\phi)}} \ge  \Pr\bb{\cc{L_1(\phi')}} \ge (1-p)\Pr\bb{\cc{L_1(\phi)}}.
\end{equation}
Similarly, since $L_2(\phi)=L_2(\phi')\cup L_x$, we have
\begin{equation}\label{zz2}
 \Pr\bb{\cc{L_2(\phi')}} \ge  \Pr\bb{\cc{L_2(\phi)}} \ge (1-p)\Pr\bb{\cc{L_2(\phi')}}.
\end{equation}

As $\phi$ is increased from $0$ to $2\pi/k$, the probability
of $\cc{L_1(\phi)}$ decreases, while the probability of $\cc{L_2(\phi)}$ increases.
These probabilities change only at a finite set of jumps, corresponding to boundary
points of $D_R$ on edges of $G$, and we have just shown (in \eqref{zz1} and \eqref{zz2})
that, at each jump, $\Pr\bb{\cc{L_2(\phi)}}$ increases by at most a factor
of $(1-p)^{-1}$, and $\Pr\bb{\cc{L_1(\phi)}}$ decreases by at most this factor.
When $\phi=0$ the arc $A_1$ is empty, so $1=\Pr\bb{\cc{L_1(\phi)}} \ge \Pr\bb{\cc{L_2(\phi)}}$,
while for $\phi=2\pi /k$ the arc $A_2$ is empty and the inequality is reversed.
Increasing $\phi$ gradually from $0$ until $\Pr\bb{\cc{L_2(\phi)}}$
first exceeds $\Pr\bb{\cc{L_1(\phi)}}$, it follows that
there is some $0<\phi<2\pi/k$ for which
\begin{equation}\label{close}
 \Pr\bb{\cc{L_1(\phi)}} \le \Pr\bb{\cc{L_2(\phi)}} \le (1-p)^{-2}\Pr\bb{\cc{L_1(\phi)}}.
\end{equation}
From now on we fix such a $\phi$, and write $L_i$ for $L_i(\phi)$, and so on.

Inequalities \eqref{close} and \eqref{l12} imply
\[
 \Pr\bb{\cc{L_1}} \le \sqrt{ \Pr\bb{\cc{L_1}} \Pr\bb{\cc{L_2}} } \le 
 \eps^{1/2k}
\]
and
\[
 \Pr\bb{\cc{L_2}} \le \sqrt{ (1-p)^{-2}\Pr\bb{\cc{L_1}} \Pr\bb{\cc{L_2}} } \le 
 (1-p)^{-1}\eps^{1/2k}.
\]
Since $\Pr\bb{L_{i+2}}=\Pr\bb{L_i}$ by symmetry, it follows that
for every $i$ we have
\[
 \Pr\bb{\cc{L_i}} \le (1-p)^{-1}\eps^{1/2k} \le 1/5,
\]
where the final inequality is from \eqref{edef}.
In other words, $\Pr(L_i)\ge 4/5$ for every $i$.

From \eqref{l12s}, there is {\em some} $j\in \{1,2\}$ for which
$\Pr\bb{\cc{L_j^\star}}\le \eps^{1/2k}\le 1/5$.
By symmetry, $\Pr(L_{j+2}^\star)=\Pr(L_j^\star)\ge 4/5$.
Let $E=L_j^\star\cap L_{j+1}\cap L_{j+2}^\star\cap L_{j+3}$.
Then $E$ is the intersection of four events each of which has probability
at least $4/5$, so $\Pr(E)\ge 1/5>0$.

To complete the proof of Theorem~\ref{th1}, we shall invoke the
uniqueness result of Aizenman, Kesten and Newman~\cite{AKNunique} mentioned
in the introduction. This result states that, if $\Lambda$
is a connected $d$-dimensional lattice, and each edge $e$ of $\Lambda$ is open with probability
$0<p_e<1$ independently of the other edges, then, provided
the edge probabilities $p_e$ are preserved by translations
of the lattice into itself,
the probability that there are two or more infinite open clusters is zero.
These assumptions apply both to $G$ and to $\Gs$, so, 
with probability $1$, neither $G$ nor $\Gs$ contains two infinite open
clusters. Thus, with positive probability, $E$ holds and the infinite
open paths $P_{j+1}$, $P_{j+3}$ in $G$ leaving $D_R$ from $A_{j+1}$
and $A_{j+3}$ are joined by a finite open path $P$ in $G$,
while the infinite open paths $P_j^\star$, $P_{j+2}^\star$ in $\Gs$
leaving $D_R$ from $A_j$ and $A_{j+2}$ are joined by a finite open path $P^\star$
in $\Gs$.
Since open paths in $G$ and in $\Gs$ cannot cross, this is a topological impossibility.
This contradiction completes the proof.
\end{proof}

\section{Consequences and extensions}

In this section we give a corollary of Theorem~\ref{th1}, and also
some extensions. In all these, the reduced symmetry requirement of Theorem~\ref{th1}
itself is the only new ingredient, so we shall only outline the arguments.

Menshikov~\cite{Menshikov} proved (essentially) that if $G$ is a lattice,
and we take the bonds of $G$ to be open independently with probability $p<\pHb(G)$,
then as the graph distance between two vertices of $G$ increases, the probability
that they are joined by an open path decreases exponentially.
(This also follows from an independent result of Aizenman and Barsky~\cite{AiBar87};
see~\cite{GG99} or~\cite{BRbook}.)
A well known consequence of this result is
that if $G$ is a plane lattice and $\Gs$ is its dual, then $\pHb(G)+\pHb(\Gs)\le 1$.
(See~\cite{BRbook}, for example.) 
Theorem~\ref{th1} implies that $\pHb(G)+\pHb(\Gs)\ge 1$,
so we obtain the following corollary.

\begin{corollary}\label{c1}
Let $G$ be a plane lattice with $k$-fold symmetry, $k\ge 2$, 
and let $\Gs$ be its dual. Then $\pHb(G)+\pHb(\Gs)=1$.
\hfill $\square$
\end{corollary}

Theorem~\ref{th1} and Corollary~\ref{c1} have natural equivalents
for site percolation. In {\em independent site percolation} on a graph
$G$, we take each vertex of $G$ to be open with probability $p$,
independently of the other vertices. Open paths and clusters
are then paths and components in the subgraph of $G$ induced by the open
vertices, and the definitions of the percolation probability $\theta_G^\site(p)$
and critical probability $\pHs(G)$ are analogous to those of $\theta_G^\bond(p)$
and $\pHb(G)$.

There is a well known `matching' relation
between certain pairs of graphs that plays a role for site percolation corresponding
to plane duality for bond percolation, used by
Kesten~\cite{KestenBook}, for example, but going back much earlier.
A plane lattice $G$ matches the (not necessarily planar) graph $\Gd$ obtained from $G$
by adding an edge $xy$ for every pair of non-adjacent vertices $x$, $y$
of $G$ lying in a common face of $G$. The relevant properties of $\Gd$
are that finite open clusters in $G$ or $\Gd$ are surrounded by closed
cycles in $\Gd$ or $G$, respectively, and that an open path in $G$
cannot cross a closed path in $\Gd$. 

Apart from the minor irritation that all references to open paths or
clusters in the dual graph (now $\Gd$ rather than $\Gs$) should be
changed to closed paths or clusters, the proofs of Theorem~\ref{th1}
and Corollary~\ref{c1}
carry over almost verbatim to the site percolation setting, to give
the following result.

\begin{theorem}\label{th2}
Let $G$ be a plane lattice with $k$-fold symmetry for some $k\ge 2$,
and let $0<p<1$. Then either $\theta_G^\site(p)=0$ or $\theta_{\Gd}^\site(1-p)=0$.
Furthermore, $\pHs(G)+\pHs(\Gd)=1$.
\hfill$\Box$
\end{theorem}

Just like Zhang's original argument, the proof of Theorems~\ref{th1} and~\ref{th2}
carries over to weighted graphs, in which edges are open independently
but different edges have different weights, or probabilities of being open.
Of course, these probabilities should respect the lattice structure. For a
formal statement and proof see \cite[Section 5.4]{BRbook}.
In this way one can obtain a simple proof of a result of Grimmett~\cite{GG99}
(stated, but not proved, by Kesten~\cite{KestenBook})
giving the `critical surfaces'
for the weighted triangular and hexagonal lattices, in which
the weight of an edge depends on its orientation; see \cite[p.~152]{BRbook}.

Like Zhang's argument, the proof of Theorem~\ref{th1} (which,
after all, is just a simple modification of Zhang's argument) can be applied
in contexts other than ordinary percolation. In particular, because the proof
relies only on the uniqueness of the infinite cluster, positive correlation of increasing
events, and the probability of a single edge being open being bounded away from $1$,
it also applies to the random cluster model of Fortuin and Kasteleyn~\cite{FK}.
Indeed, Welsh~\cite{W93}
showed that Zhang's argument can be used to prove
that the critical probability $\pH(q)$ for the random cluster model on the square lattice
satisfies
\begin{equation}\label{ZRC}
 \pH(q)\ge \sqrt{q}/(1+\sqrt{q})
\end{equation}
for all $q\ge 1$ (see also Grimmett~\cite{GG95,GGRC});
as far as we are aware, this is the only known proof of this result.
The proof of Theorem~\ref{th1} extends this result to the random cluster
model on any self-dual plane lattice with $k$-fold symmetry. More generally,
it shows that percolation cannot occur simultaneously in two `dual' random cluster measures
with $k$-fold symmetry.

The latter observation answers a conjecture of Grimmett~\cite[relation (6.73)]{GGRC}.
For ${\bf p}=(p_1,p_2,p_3)$ and $q\ge 1$, let $\phi^0_{T,{\bf p},q}$ denote the 
(lower) anisotropic random cluster measure on the triangular lattice, in which edges
receive weights $p_1$, $p_2$ or $p_3$ according to their direction, and a weight
$q$ is assigned to each cluster; see~\cite{GGRC} for a formal definition, which takes
a little time to set up. Also, let $\theta^0_T({\bf p},q)$ denote the probability
that the origin is in an infinite open cluster in this measure.
\begin{theorem}\label{thG}
Let $q\ge 1$ and $0<p_1,p_2,p_3<1$ be such that
\begin{equation}\label{RCdual}
 y_1y_2y_3+y_1y_2+y_2y_3+y_3y_1-q \le 0,
\end{equation}
where $y_i=p_i/(1-p_i)$.
Then $\theta^0_T({\bf p},q)=0$.
\end{theorem}
\begin{proof}
As noted by Grimmett~\cite{GGRC}, the measure $\phi^0_{T,{\bf p},q}$ is dual
to a corresponding measure $\phi^1_{H,{\bf r},q}$ on the hexagonal lattice, with weights $r_i$
satisfying $r_i/(1-r_i)=q(1-p_i)/p_i$.
Grimmett also shows that if equality holds in~\eqref{RCdual}, then,
using the star-triangle inequality,
this latter measure is equivalent to a measure on the triangular lattice that stochastically
dominates $\phi^0_{T,{\bf p},q}$. We may assume
equality in \eqref{RCdual} by increasing $p_1$, say; the resulting random cluster measure
stochastically dominates the original.
It follows that if $\theta^0_T({\bf p},q)>0$, then
we have percolation in two dual random-cluster measures on lattices
with $3$-fold symmetry, contradicting the analogue of
Theorem~\ref{th1}, whose proof readily adapts to this setting.
\end{proof}

%Using our modification of Zhang's argument, the proof of~\eqref{ZRC}
%or of its analogue for the inhomogeneous square lattice~\cite[Theorem 6.54]{GGRC}
%extends immediately to prove this conjecture.

\medskip
Let us give one further application of the proof of Theorem~\ref{th1}, proving a result
we believe is new.
Let $0\le r\le 1$ be fixed, and let $T_r$ be the random triangulation of the plane
obtained from the square lattice $\Z^2$ as follows: add one diagonal to every face,
choosing the diagonal parallel to $(1,1)$ with probability $r$, and that parallel
to $(1,-1)$ with probability $1-r$, these choices being independent for different faces.
If $r=0$ or $r=1$, then $T_r$ is isomorphic to the triangular lattice $T$,
which satisfies $\pHs(T)=1/2$.
\begin{theorem}\label{th4}
Suppose that $0\le r\le 1$. With probability $1$, the random triangulation
$T_r$ satisfies $\pHs(T_r)=1/2$.
\end{theorem}
\begin{proof}
For $0\le r\le 1$ and $0<p<1$, first construct $T_r$ as above,
and then select vertices of $T_r$, i.e., points of $\Z^2$,
to be open with probability $p$, independently of each other and of $T_r$.
Let $\Pr_{r,p}$ denote the associated probability measure, and $\theta(r,p)$
the probability that the origin is in an infinite open path.
Our aim is to adapt the proof of Theorem~\ref{th1}, or, rather, of
the site percolation version, Theorem~\ref{th2}, to show that
\begin{equation}\label{cl}
 \theta(r,1/2)=0
\end{equation}
for any $r$.

Any realization of the random graph $T_r$ is a triangulation,
and any triangulation $G$ satisfies $\Gd=G$, so open paths in the random
graph $T_r$ cannot cross closed paths in $T_r$, and so on. Furthermore,
when $p=1/2$ the distribution of the open subgraph of $T_r$ is identical
to that of the closed subgraph. Since the model associated to $\Pr_{r,p}$
has $2$-fold symmetry, the proof of Theorem~\ref{th2} would apply, {\em mutatis
mutandis}, to prove \eqref{cl}, except for one problem: we do not
have the required positive correlation. For example,
the existence of the single-edge paths $(0,0)(1,1)$ and $(0,1)(1,0)$
are negatively correlated (indeed, mutually exclusive) events.
It turns out that we can get around this.

Let $X$ be the graph obtained
from $\Z^2$ by adding two new vertices $v_{F,1}$, $v_{F,2}$
for each face $F$ of $\Z^2$. Join $v_{F,1}$ to the top-right and bottom-left
vertices of $F$, and $v_{F,2}$ to the other two vertices. Let $v_{F,1}$
be open with probability $r$, and $v_{F,2}$ with probability $1-r$,
and each vertex of $\Z^2$ be open with probability $p$, with
all these events independent.
When considering events depending only on which vertices of $\Z^2$ are
joined by open paths, the
(non-planar) lattice site percolation measure just defined is
equivalent to $\Pr_{r,p}$. Indeed, having examined the states
of all vertices of $\Z^2$, we need only look within a face $F$
(i.e., consider the vertices $v_{F_,i}$ in $X$ or consider the diagonals
in $T_r$) if two opposite vertices of $F$ are open, and the other two
vertices closed. (Otherwise, all open vertices of $F$ are joined by
paths using edges of $\Z^2$.) In this case, the conditional probability
that the relevant vertex of $X$ is open is exactly the conditional probability
that the relevant diagonal of $T_r$ is present.

Applying Harris' Lemma to site percolation on $X$, it follows that
events such as `there is an open path in $T_r$ leaving the disc $D_R$ from a certain arc'
are positively correlated. As noted above, the proof of Theorem~\ref{th2} thus goes through,
to show that the  $\Pr_{r,1/2}$-probability that the origin is in an infinite path
in $T_r$ is 0. Hence, $\pHs(T_r)\ge 1/2$ with probability $1$.
Applying a suitable form of Menshikov's Theorem to the lattice $X$,
it is easy to show that $\pHs(T_r)\le 1/2$ with probability 1, and Theorem~\ref{th4} follows.
\end{proof}
\medskip
We close with a question: can one adapt Zhang's argument, or the proof of Theorem~\ref{th1},
to the case of a lattice with mirror symmetry, but without $k$-fold symmetry for any $k$?
Past experience suggests that it would be dangerous to conclude that no such variant of the
proof exists purely on the basis that we have failed to find one!
As noted earlier, Sheffield~\cite{Sheffield} has proved such a result requiring
only the symmetry of a lattice; his proof is far from simple, however.

\medskip\noindent
{\bf Acknowledgements.}
The writing of the first draft of
this paper was completed during a visit of
the authors to the Institute for Mathematical Sciences,
National University of Singapore. We are grateful to the Institute
for its hospitality.
We should also like to thank Geoffrey Grimmett for drawing our attention
to his conjecture answered by Theorem~\ref{thG}.

\end{document}